\documentclass[11pt]{amsart}
\usepackage[margin=1.0in]{geometry}
\usepackage{amsmath,amssymb,amsxtra,amsthm,hyperref,color}
\usepackage{tikz}



\newcommand{\bh}[1] {\mathcal{B}(\mathcal{#1})}
\newcommand{\pref}[1] {(\ref{#1})}

\newcommand{\lspan} {\operatorname{span}}

\newcommand{\cH}{\mathcal{H}}
\newcommand{\cK}{\mathcal{K}}

\newtheorem{theorem}{Theorem}[section]
\newtheorem{corollary}[theorem]{Corollary}
\newtheorem{proposition}[theorem]{Proposition}    
\newtheorem{lemma}[theorem]{Lemma}

\theoremstyle{definition}
\newtheorem{definition}[theorem]{Definition}

\newtheorem{example}[theorem]{Example}
\newtheorem{remark}[theorem]{Remark}

\numberwithin{equation}{section}
\def\N{\mathbb N}
\def\C{\mathbb C}

\def\Z{\mathbb Z}

\begin{document}

\title[Amenability and functoriality of semigroup C*-algebras]{Amenability and functoriality  of  \\
right-LCM semigroup C*-algebras}

\author{Marcelo Laca}
\address{Department of Mathematics and Statistics, University of Victoria, Victoria, B.C. V8W 3R4}
\email{laca@uvic.ca}
\author{Boyu Li}
\address{Department of Mathematics and Statistics, University of Victoria, Victoria, B.C. V8W 3R4}
\email{boyuli@uvic.ca}
\date{\today}


\begin{abstract} We prove a functoriality result for the full C*-algebras of right-LCM monoids 
with respect to monoid inclusions that are closed under factorization and preserve orthogonality, 
and use this to show that if a right-LCM monoid is amenable in the sense of Nica, then so are its submonoids. 
As applications, we complete the classification of Artin monoids with respect to Nica amenability by showing that
only the right-angled ones are amenable in the sense of Nica and we show that the Nica amenability of a graph product 
of right-LCM semigroups is inherited by the factors.
\end{abstract}

\maketitle

\section{Introduction}
The study of C*-algebras generated by semigroups of isometries goes back over half a century, and has its origins in the seminal papers of L. Coburn on the C*-algebra of the unilateral shift \cite{Coburn1967},  of R. Douglas on C*-algebras of discrete subsemigroups of the positive reals \cite{Douglas1972}, and of J.~Cuntz on C*-algebras generated by isometries with mutually orthogonal ranges \cite{Cuntz1977}. Their work revealed early on that the theory is deeply rooted in noncommutative function theory and noncommmutative dynamical systems; a case in point being the role of von Neumann's inequality in Coburn's proof of uniqueness of the C*-algebra generated by a nonunitary isometry; another one being the role of the shift on a UHF algebra in Cuntz's proof of uniqueness of $\mathcal O_n$. 

As with groups, the primary interest lies in the C*-algebra generated by the left regular representation of a semigroup on its Hilbert space, but there is also a very useful C*-algebra that is universal with respect to representations by isometries that satisfy a certain covariance relation abstracted from the left regular representation. Much of the research has been driven by the uniqueness question, of whether the algebraic relations in the semigroup determine the structure of its C*-algebra, and, in particular, the related question of whether the left regular representation has the universal property itself. 

In his study of quasi-lattice ordered groups \cite{Nica1992}, A. Nica realized that these two key questions of uniqueness and faithfulness depend on an amenability property of the semigroup and on a properness condition of the representation; the theme was subsequently picked up in \cite{LacaRaeburn1996}, where it was formalized into a method for establishing amenability and for verifying properness by casting the semigroup C*-algebras as semigroup crossed products. 
In rough terms, the method to study amenability consists of exhibiting a factorization of the canonical conditional expectation of a semigroup crossed product
that shows it is faithful as a positive map, which is equivalent to Nica amenability.
 It is a bit surprising that many nonamenable groups have natural subsemigroups that are in fact Nica amenable; for instance, as already observed by Nica, the free semigroup has this property, which manifests itself as the uniqueness of the C*-algebra generated by $n$ isometries with orthogonal ranges that are jointly proper in the sense that the orthogonal complements of their ranges have nontrivial intersection, a result proved by Cuntz in \cite{Cuntz1981}.

There has been a lot of activity on C*-algebras generated by semigroups of isometries during the last two decades, partly because they produce 
tractable yet interesting families of examples. Arguably, the main achievement has been the very successful theory developed by Xin Li for general left-cancellative semigroups by way of a covariance condition expressed in terms of constructible ideals, see  \cite{XLi2012,CELY2017} and references thereof. 
Here we are interested in the particular case of semigroups for which the constructible ideals are all principal right ideals and thus are
easily accessible by way of a slight generalization of Nica's ideas \cite{Starling2015, BLS2017}, so we will not need the full power of Xin Li's theory.
We now know right-angled Artin monoids are Nica amenable, but nonabelian Artin monoids of spherical type (also known as finite-type) are not \cite{CrispLaca2002};  the ideal structure and the K-theory of the C*-algebras and the boundary quotients for right-angled Artin monoids have been computed \cite{CrispLaca2007,Ivanov2010,ELR2016};  new classes of examples including the Baumslag-Solitar semigroups have been proven to be
quasi-lattice ordered and Nica amenable \cite{HRT2018}, and more recently, A. an Huef, B. Nucinkis, C.F. Sehnem, and D. Yang have 
used a generalization of the controlled maps from  \cite{LacaRaeburn1996} and \cite{HRT2018} to produce more Nica amenable semigroups 
and have given criteria for nuclearity of the semigroup C*-algebras \cite{ABCD2019}.

The two classes of Artin monoids mentioned above, namely right-angled (Nica amenable) and non-abelian spherical type (not Nica amenable) are in a sense opposite extremes in the wide spectrum of Artin monoids, and the natural question of Nica amenability for the intermediate types of Artin monoids
was raised at the end of \cite{CrispLaca2002}. The original motivation for the present work is to answer this question, and thus complete the classification of Artin monoids with respect to Nica amenability.  A word is in order to explain the context.
We  know by now that all Artin monoids embed in their groups  \cite{Paris2002}, 
but we only know that those embeddings  give rise to quasi-lattice ordered groups   
in the spherical and right-angled cases. However, every Artin monoid satisfies the 
right LCM property by the reduction lemma of Brieskorn and Saito \cite[Lemma 2.1]{BrieskornSaito1972} and since a quasi-lattice ordered group structure is not 
really necessary to consider covariant isometric representations in general, and the left regular representation in particular, 
the concept of Nica amenability makes sense for right LCM monoids.    
Thus, we begin with a preliminary section in which we collect and update the relevant definitions and characterizations of covariance and amenability in the sense of Nica for C*-algebras of left-cancellative right LCM monoids. 
 
Since Nica amenability is a property of C*-algebras constructed from semigroups, it is not surprising that functoriality of the construction is at the heart of the question. 
Functoriality results specific to right-LCM semigroups have been obtained in \cite{BLS2018}, and include criteria to decide whether a homomorphism of right-LCM semigroups lifts to a C*-algebra homomorphism of their universal C*-algebras for covariant representations. The difficulty here is that, in contrast to what happens for groups, it is not clear whether
the C*-algebra homomorphism resulting from the inclusion of semigroups is injective. In Section \ref{sec:main} we establish this injectivity for inclusions of right-LCM semigroups that are {\em closed under factorization and preserve orthogonality}, see Definitions  \ref{def:closeduf} and \ref{def:preservesorth} for the precise formulation of these concepts. Through a new direct approach that relies on recent developments in dilation theory for right-LCM semigroups \cite{BLi2019}, we then prove that a covariant representation of a submonoid dilates/extends to a covariant representation of the larger monoid on a larger Hilbert space, Theorem \ref{thm:diln}, 
which is our main technical result. 
We  then establish in Corollary~\ref{cor.functor.univ} that, under the suitable factorization and orthogonality assumptions, an inclusion of right-LCM monoids gives an injective homomorphism of the full semigroup $C^*$-algebra.  This, in turn, implies that Nica amenability is inherited by  subsemigroups, Theorem~\ref{thm.amenable}. 

In Section \ref{sec:applications} we apply our main results to show that, with the exception of right-angled Artin monoids, which are known to be Nica amenable by \cite[Theorem 24]{CrispLaca2002}, all the remaining Artin monoids are not Nica amenable, Theorem \ref{thm.Artin}, so  their full and  reduced C*-algebras are effectively different. This means that two different covariant isometric representations of an Artin monoid that is not right-angled do not necessarily determine isomorphic C*-algebras even if they are jointly proper.
As a second application of our results we also show that a graph product of right-LCM monoids is Nica amenable only if all the factors are Nica amenable Corollary~\ref{cor:graphprod}, which is a sort of converse, for cancellative right LCM semigroups, of a known criterion for amenability, cf. \cite[Theorem 20]{CrispLaca2002} and \cite[Corollary 5.27]{ABCD2019}.
 
 The authors would like to thank Nadia Larsen for pointing out the relevance of \cite{BLS2017} and Xin Li for several helpful comments. 
 
\section{Preliminaries on semigroups and isometric representations}

A semigroup $P$ is a set with an associative multiplication.
We assume that $P$ is a {\em monoid}, which means that $P$ contains an element $e$ such that $e\cdot x=x\cdot e=x$ for all $x\in P$; 
we leave open the possibility for 
$P$ to also contain nontrivial {\em units}, i.e.  elements that are invertible in $P$, and we denote the set of units by $P^*$. 
We shall also assume throughout that the semigroup $P$ is  {\em left-cancellative}: this means that if $p\cdot x=p\cdot y$ for  $x,y,p\in P$, then  $x=y$; 
left cancellation ensures that the canonical left regular representation of $P$ on $\ell^2(P)$ is by isometries.
For our main results about amenability, we will need to assume that $P$ is also {\em right-cancellative}; this has the obvious analogous definition
and allows us to identify the diagonal subalgebras as the range of conditional expectations that lie at the heart of our approach.
Needless to say, both cancellation properties hold whenever $P$ embeds in a group, even though such a group is not uniquely determined and there is often no canonical choice for it. 

 A  {\em right-LCM semigroup}  is a left-cancellative semigroup $P$ such that for every pair $p,q\in P$ either $pP\cap qP =\emptyset$ or there exists an element $r\in P$ such that $pP\cap qP=rP$, in other words, the nonempty intersection of two principal right ideals is another principal right ideal (alternatively, when this happens one says that $P$ {\em satisfies Clifford's condition}, cf. \cite[Definition 2.18]{Norling2014}). Such an element $r$ is called a {\em least common multiple} of $p$ and $q$, and is not unique if and only if $P$ has nontrivial units, that is, if and only if $P^*\neq \{e\}$. Indeed, 
 if $u\in P^*$ then $rP =  ruP$ so both $r$ and $ru$ are LCMs of the same pairs of elements; conversely, if $r$ and $s$ are both least common multiples of $p$ and $q$, then $rP=sP$, in which case an easy computation using left-cancellation shows that $r=s\cdot u$ for some unit $u\in P^*$. 
 
We use $p\vee q:=\{r: pP\cap qP=rP\}$ to denote the set of all least common multiples of $p$ and $q$, and note that $p\vee q = r P^*$ for any given $r\in p\vee q$. When $F\subset P$ is a finite subset, we also use $\vee F:=\{r: rP=\bigcap_{p\in F} pP\}$ to denote the set of all least common multiples of the elements of $F$. Notice that $\vee F$ can be empty: for example, the set of generators in the free semigroup on two or more generators has no common upper bound. When $p \vee q \neq \emptyset$ for every pair $p,q \in P$ we say that $P$ is a semi-lattice. 
The possibility of $P$ being a group itself, that is $P = P^*$, is not presently ruled out, but we will have nothing new to say in this case.

An important class of right-LCM semigroups is the class of positive cones in quasi-lattice ordered groups \cite{Nica1992}, and, more generally, the positive cones in the
`weakly quasi-lattice ordered groups' recently introduced in \cite{ABCD2019}, which, in our terminology, would be described as
 {\em group-embeddable right-LCM semigroups with no nontrivial units}.
The extra generality arises from only requiring the quasi-lattice property on the semigroup $P$, and not on the whole group, see \cite[Remark 8]{CrispLaca2002}.   We refer to \cite[Section 3]{ABCD2019} for motivation and details of weakly quasi-lattice ordered groups, and for an example that does not give a quasi-lattice ordered group in any embedding, \cite[Proposition 3.10]{ABCD2019}. We point out that the semigroups considered here  are more general in that they are not assumed to embed in a group and  they may contain a nontrivial group of units.  

Since $P$ is left cancellative, for each fixed $p\in P$ the map $q\mapsto {pq}$ 
is a bijection of $P$ onto  $pP$; thus, denoting by $\{\delta_q\}_{q\in P}$ the standard orthonormal basis of $\ell^2(P)$, 
we may extend $\lambda(p): \delta_q \mapsto \delta_{pq}$ by linearity and continuity to obtain an isometry $\lambda(p)$ on $\ell^2(P)$, with adjoint given by 
\[
\lambda(p)^* \delta_q = \begin{cases}\delta_{q'} & \text{ if } q = pq' \in pP\\
0 &\text{ if } q\notin pP.
\end{cases}
\]
It is easy to see that $\lambda(p) \lambda(q) = \lambda(pq)$ and that the projection  $\lambda(p)\lambda(p)^* $ is the multiplication operator by the 
characteristic function of the set $pP$. The representation $\lambda : P \to \mathcal B(\ell^2(P))$ is called the {\em left regular representation} of $P$.
Since $\lambda(p)$  is a unitary when $p\in P^*$,  this generalizes the usual left regular representation of a group.
Inspired by the properties of the left regular representation of the positive cone in a quasi-lattice ordered group, 
Nica isolated the following notion of covariance for isometric representations \cite{Nica1992}, which also works for right LCM monoids.

\begin{definition}An isometric representation $V:P\to\bh{H}$ of a right-LCM semigroup $P$ is  {\em Nica-covariant} if for every $x,y\in P$, 
\begin{equation}\label{eqn:covariance}V_xV_x^* V_yV_y^* = \begin{cases} 
V_r V_r^* &\mbox{ if } r\in x\vee y\\
0, &\mbox{ if } x\vee y =\emptyset.
\end{cases}
\end{equation}
To see that the definition is independent of the choice of $r$, recall that any two least common upper bounds $r,s\in x\vee y$ of $x$ and $y$ in $P$ satisfy $r=su$ for some invertible $u\in P^*$, hence $V_r V_r^* = V_sV_u V_u^*V_s^* = V_sV_s^*$ because  $V_u$ is an invertible isometry, namely a unitary element.
  \end{definition}
\begin{example} If $V$ is an isometric representation of the monoid $P = \N^2$, then $V$ determined by the two commuting isometries $S:=V_{(1,0)}$ and $T:=V_{(0,1)}$, and $V$ is  Nica covariant if and only if $S^*T = TS^*$, which is satisfied by e.g. the left regular representation but not by the choice $S=T$. 
On the other hand, if  $V$ is a representation of the monoid $P =F^+_2$, with generators $a$ and $b$,
then $V$ is Nica covariant if and only if $V_a$ and $V_b$ have mutually orthogonal ranges.
\end{example} 

In a close analogy to what happens for groups, there are two C*-algebras associated with $P$: one is the {\em reduced $C^*$-algebra}  $C^*_\lambda(P)$, which, by definition,  is the C*-algebra of operators on $\mathcal B(\ell^2(P))$ generated by the left regular representation of $P$ on $\ell^2(P)$, 
and is the concrete object of interest. The other one is
the {\em full semigroup $C^*$-algebra} $C^*(P)$, which, by definition, is the universal $C^*$-algebra generated by a Nica-covariant representation 
$V^u$; and is an abstract object whose representation theory encodes the representations of $P$ by covariant semigroups of isometries.
Since the left regular representation is obviously Nica-covariant, the map $V^u_p \mapsto \lambda(p)$ for  $p\in P$ extends to  a surjective $*$-homomorphism 
$C^*(P)\to C^*_\lambda(P)$. In general, this homomorphism is not an isomorphism and we have the following definition for right-LCM monoids,  extending one of the C*-algebraic characterization of amenability of a discrete group.

\begin{definition} A right-LCM monoid $P$ is called {\em Nica amenable} if the left-regular representation $\lambda:C^*(P)\to C^*_\lambda(P)$ is an isomorphism
of C*-algebras. 
\end{definition} 
\begin{example} Let  $\mathbb{F}_k^+$ be the free monoid on $k$ generators; it is easy to verify that $(\mathbb{F}_k, \mathbb{F}_k^+)$ is a quasi-lattice ordered group, which is  Nica amenable by \cite[5.1 Proposition]{Nica1992}. 
\end{example} 

Multiplying the covariance condition \eqref{eqn:covariance}  by $V_x^*$ on the left and $V_y$ on the right gives 
\[V_x^* V_y = \begin{cases} 
V_{z_1} V_{z_2}^* &\mbox{ if } z\in x\vee y, \text{ with } z=x z_1=yz_2\\
0, &\mbox{ if } x\vee y =\emptyset,
\end{cases}\]
from which it follows that ${\lspan}\{V_pV_q^*: p,q\in P\}$ is a dense $*$-subalgebra of
the $C^*$-algebra generated by the collection $\{V_p\}_{p\in P}$.
It also follows from  \eqref{eqn:covariance}  that  $\mathcal{D}_P:=\overline{\lspan}\{V^u_p (V^u_p)^*: p\in P\}\subset C^*(P)$ is a commutative C*-subalgebra of $C^*(P)$. The canonical image of $\mathcal{D}_P$ under $\lambda$ is the C*-algebra $ \mathcal{D}_{P,\lambda}: = \overline{\lspan}\{\lambda(p) \lambda(p)^*: p\in P\} \subset C_\lambda^*(P)$, which can be regarded as the subalgebra of $\ell^\infty(P)$ generated by the characteristic functions of the sets $pP$, acting by multiplication on $\ell^2(P)$.
 We will focus on situations where the diagonal algebras coincide and are the range of conditional expectations; 
Nica amenability can then be characterized as faithfulness of the conditional expectation from $C^*(P)$ onto $\mathcal{D}_P$.
 The following proposition summarizes and combines various results from
 \cite{XLi2012,Nica1992,Norling2014}, and especially from \cite[Proposition 3.14]{BLS2017}; we include it here with its short proof mainly to indicate
 that the various arguments do not require $P$ to embed in a  group and work for right-LCM monoids that are cancellative on both sides.
\begin{proposition} \label{pro:conditionalexpectations}
Suppose P is a (right-left) cancellative right-LCM monoid. Then 
\begin{enumerate}
\item the map $\Phi_\lambda: \lambda(p) \lambda(q)^* \mapsto \delta_{p,q}\lambda(p) \lambda(p)^*$ for $p,q \in P$ extends to a surjective conditional expectation 
$\Phi_\lambda:C^*_\lambda(P)\to \mathcal{D}_{P,\lambda}$ that is always faithful on positive elements;
\item the map  $\Phi: V^u_p (V^u_q)^* \mapsto \delta_{p,q}V^u_p (V^u_p)^*$  for $p,q \in P$ extends to a surjective conditional expectation $\Phi:C^*(P)\to \mathcal{D}_{P}$
\item $\Phi$ is faithful on positive elements if and only if   $P$ is Nica amenable.
\end{enumerate}
\end{proposition} 
\begin{proof}
It is  well known that if $\mathcal H$ is a Hilbert space with orthonormal basis $\{\delta_z\}_{z\in I}$ and if we denote by $E_z$ the rank-one projection onto $\C\delta_z$,
then the diagonal map defined by $\Delta (T) := \sum _z E_z T E_z$ for $T\in \bh{H}$ is a unital, completely positive, idempotent, linear, contractive map onto the C*-algebra of bounded linear operators on $\mathcal H$ that have a diagonal matrix in the given basis; in addition, $\Delta$ is faithful as a 
conditional expectation, see e.g. the Remark in \cite[3.6]{Nica1992}. So in order to prove part (1)  it suffices to verify  that  $\Delta(\lambda(p) \lambda(q)^*)  = \delta_{p,q}\lambda(p) \lambda(p)^*$ for all $p, q \in P$.
 Evaluation of the left hand side at an orthonormal basis element $\delta_a$ gives
\begin{align}
\Delta(\lambda(p) \lambda(q)^*) \delta_a &= 
E_a \lambda(p) \lambda(q)^*E_a \delta_a\\
 &= \begin{cases} E_a \delta_{pa'} & \text{ if } a = qa'\in qP\\
0 & \text{ if } a\notin qP
\end{cases}\\
&= \begin{cases} \delta_a  & \text{ if } a = qa'\in qP \text{ and } pa' = qa'\\
0 & \text{ if } a\notin qP \text{ or } a =qa'\in qP \text{ but } pa' \neq qa',
\end{cases}
\end{align}
where in the case $a\in qP$ we have denoted by $a'$ the unique element in $P$ such that $a = qa'$.
Evaluation of the right hand side at $\delta_a$ gives
\[
\delta_{p,q}\lambda(p) \lambda(q)^* \delta_a =  \begin{cases}  \delta_a & \text{ if } a = qa'\in qP \text{ and } p = q  \\
0 & \text{ if } a\notin qP  \text{ or } a = qa' \in qP \text{ but }p\neq q.
\end{cases}
\]
Since $P$ is assumed to be right cancellative, $pa' = qa'$ is equivalent to $p = q$ and the cases match, so we conclude that the restriction of $\Delta$ to 
$C^*_\lambda(P)$ is the required map $\Phi_\lambda$; this is the only part of the argument that uses right cancellation. 
It is clear that the range of $\Phi_\lambda$ is the closed linear span of the projections $\lambda(p)\lambda(p)^*$, 
namely the diagonal $\mathcal D_{P,\lambda}$.

In order to prove part (2) we use a neat trick from the proof of Lemma 2 and Definition in \cite[4.3]{Nica1992}.
Recall that left cancellative right-LCM monoids satisfy the independence condition by 
 \cite[Proposition 2.19]{Norling2014}, so we know from \cite[Corollary 2.27]{XLi2012} that 
 the restriction of the left regular representation to the diagonal algebra gives an isomorphism $\mathcal D_P \cong \mathcal D_{P,\lambda}$,
 see also \cite[Lemma 3.18]{BLS2017} for semigroups without an identity element. 
Hence, following Nica, we define $\Phi$ by going around the long way and reversing the bottom horizontal arrow in the diagram 
\begin{equation*}
    \centering
    \begin{tikzpicture}[scale=0.9]
    \node at (-2,2) {$C^*(P)$};
    \node at (2,2) {$C^*_\lambda(P)$};
    \node at (-2,0) {$\mathcal{D}_P$};
    \node at (2,0) {$\mathcal{D}_{P,\lambda}$};
    
    \draw[->] (-1,2) -- (1,2);
    \draw[->] (-1,0) -- (1,0);
	\draw[->] (-2,1.5) -- (-2,0.5);
    \draw[->] (2,1.5) -- (2,0.5);
    
    \node at (-2.5, 1) {$\Phi$};
    \node at (2.5, 1) {$\Phi_\lambda$};
    \node at (0, 2.5) {$\lambda$};
    \node at (0, 0.5) {$\lambda|_{\mathcal{D}_P}$};
    \end{tikzpicture}
\end{equation*}
which ensures at one stroke that $\Phi$ is a conditional expectation mapping $V^u_p (V^u_q)^*$ to $\delta_{p,q}V^u_p (V^u_p)^*$ 
and that the diagram commutes.
A standard argument then shows that Nica amenability, that is to say,  faithfulness of the left regular representation $\lambda$ of $C^*(P)$ is equivalent to faithfulness of $\Phi$ on positive elements, proving part (3).
\end{proof}
Nica amenability has been established for  several classes of right-LCM semigroups, and the most common strategy has been
to use generalizations of the notion of controlled
map originally introduced in  \cite{LacaRaeburn1996}, see, e.g. \cite{CrispLaca2002,CrispLaca2007,HRT2018} and especially \cite{ABCD2019} for the latest version.
It is also known that Nica amenability follows from certain approximation property, \cite{Nica1992}, \cite{Exel1997}. 

If a quasi-lattice ordered group $(G,P)$ is semi-lattice ordered in the sense that the intersection of any two principal right ideals is a (nonempty) principal right ideal, then
 $P$ is Nica amenable if and only if the group itself is amenable, see \cite[Section 4.5]{Nica1992} or \cite[Lemma 6.5]{LacaRaeburn1996} for one direction and  \cite[Proposition 28]{CrispLaca2002} for the converse (note that semi-lattice orders were called lattice orders in \cite{CrispLaca2002}).
 The underlying reason for this is that representations of semi-lattices by unitary isometries
 satisfy the covariance condition,  and thus the group C*-algebra $C^*(G)$ is a quotient of $C^*(P)$. 
 In particular, this implies that all nonabelian Artin monoids of spherical type are not Nica amenable \cite[Theorem 30]{CrispLaca2002}.
 Oddly enough, with the possible exception of a few isolated examples, these semi-lattices were the only monoids known not to be Nica amenable until now. 
 
Finally, it should be noted that there are other definitions of semigroup amenability, e.g.  a semigroup $P$ is said to be left-amenable if there is a left-invariant mean on $\ell^\infty(P)$ \cite{Day1957}. This has been explored by Xin Li in the context of semigroup C*-algebras; a cancellative semigroup $P$ that satisfies the independence property is left-amenable if and only if there is a character defined on $C_\lambda^*(P)$ \cite[Theorem 5.6.42]{CELY2017}; since free semigroups are Nica amenable but do not admit characters, left-amenability is strictly stronger than Nica amenability.

\section{Main Results}\label{sec:main}
\begin{definition}\label{def:closeduf} Let $P_1$ be a submonoid of the monoid $P$.
We say that $P_1$ is {\em closed under factorization} in $P$ if  for every pair $x,y\in P$ such that $xy  \in P_1$ we have $x,y\in P_1$. 
\end{definition}
\begin{lemma} \label{lem:rlcmifcloseduf} Suppose $P$ is a right-LCM monoid and  $P_1$ is a submonoid of $P$ that is closed under factorization. Then $P_1$ is right-LCM.
\end{lemma}
\begin{proof}
If $x, y\in P_1$ and $xP_1 \cap yP_1 \neq \emptyset$, then  every $w \in xP_1 \cap yP_1$ is necessarily in $xP \cap yP = zP$  for any $z\in x\vee y$. Hence $ w = z\alpha$ for $z, \alpha \in P$ and the factorization property implies $z, \alpha \in P_1$, so $xP_1 \cap yP_1 = zP_1$.
\end{proof}
The lemma will not quite suffice for our purposes because it is conceivable that a pair of elements $x$ and $y$ in $P_1$ could have all its common upper bounds in $P\setminus P_1$, so that 
$xP_1 \cap yP_1 = \emptyset$ but $xP \cap yP \neq \emptyset$. In order to rule out this possibility, we make the following definition generalizing a notion introduced in \cite[Definition 1.1]{Crisp1999} for the case of Artin monoids.

\begin{definition}\label{def:preservesorth} Let  $P_1$ be a right-LCM submonoid of the  right-LCM monoid $P$. We say the inclusion $P_1\subseteq P$ {\em preserves  orthogonality} if  for all $x,y\in P_1$,  $xP_1\cap yP_1=\emptyset$ is equivalent to $xP\cap yP=\emptyset$. 
\end{definition} 

The following related condition appears  in \cite[Theorem 3.3]{BLS2018}, where it is used to decide whether
the inclusion $P_1 \subset P$ induces a canonical $*$-homomorphism $\varphi: C^*(P_1)\to C^*(P)$.
\begin{definition} Let  $P_1$ be a right-LCM submonoid of the  right-LCM monoid $P$. We say $P_1$ {\em respects the LCM} of $P$ if \[ xP\cap yP = (xP_1\cap yP_1) P, \qquad \text{for all } x,y\in P_1.
\] 
\end{definition}   
 One can easily verify that if  $P_1$ respects the LCM of $P$, then it must preserve orthogonality. We do not know whether the converse is true 
but the next lemma shows that the two notions do coincide for submonoids that are closed under factorization.
 Notice that the monoid $P_1=\mathbb{N}$ respects the LCM of $P=\mathbb{R}^+$, but is not closed under factorization. 

\begin{lemma}\label{lm.respectLCM} Let  $P_1$ be a submonoid of the right LCM monoid $P$
and assume that the inclusion $P_1\subset P$ is closed under factorization and preserves  orthogonality. Then $P_1$ respects the LCM of $P$. 
\end{lemma} 

\begin{proof} For any $x,y\in P_1$. Since $P_1$ is a right-LCM monoid by Lemma \ref{lem:rlcmifcloseduf}, either $xP_1\cap yP_1=\emptyset$ or $xP_1\cap yP_1=zP_1$ for some $z\in P_1$. In the first case, $xP_1\cap yP_1$ implies $xP\cap yP=\emptyset$ since $P_1$ preserves  orthogonality, in which case, 
\[ xP\cap yP = (xP_1\cap yP_1) P = \emptyset.\]
Suppose now $xP_1\cap yP_1=zP_1$, we need to verify that
\[ xP\cap yP = (xP_1\cap yP_1) P = zP.\] 
Since $P_1\subset P$ preserves orthogonality, we know 
$xP\cap yP \neq \emptyset$, so  $xP\cap yP=wP$ for some $w\in P$ by the LCM property of $P$. Since $z\in xP\cap yP = wP$, we can factor $z=wu$ for some $u\in P$. However, $z\in P_1$ and $P_1$  is closed under factorization so $w\in P_1$. Since $w\in xP$, we can factor $w=xx'$ and similarly $w=yy'$. Again,  $P_1$ is closed under factorization, so $x', y'\in P_1$ and thus $w\in xP_1\cap yP_1= zP_1$. Therefore, $wP_1=zP_1$ and $wP=zP$.   
\end{proof} 



Throughout the rest of this section, we fix a  submonoid $P_1$ of $P$ and assume that the inclusion is closed under factorization and preserves orthogonality. It then follows from Lemma \ref{lm.respectLCM} and \cite[Theorem 3.3]{BLS2018} that the inclusion $P_1 \subset P$ gives a canonical $*$-homomorphism $\varphi: C^*(P_1)\to C^*(P)$. It  is established in  that 
The map $\varphi$ is known to be isometric when both $P_1$ and $P$ are Nica amenable \cite[Proposition 3.6]{BLS2018}. We will see below that 
under the current factorization and  orthogonality assumptions,
$\varphi$ is isometric without any amenability assumptions.
The proof is an application of the following dilation theorem, which we quote here for convenience.

\begin{theorem}[Theorem 3.9 \cite{BLi2019}] \label{thm.dilation}
Let $T:P\to\bh{H}$ be a unital contractive representation of a right-LCM monoid $P$. The following are equivalent:
\begin{enumerate}
\item\label{thm.dilation.1} $T$ has a $\ast$-regular dilation;
\item\label{thm.dilation.2} $T$ has a minimal Nica-covariant isometric  dilation. In other words, there exists an isometric Nica-covariant representation $W:P\to\bh{K}$ for some Hilbert space $\mathcal{K}\supset\mathcal{H}$, so that $\mathcal{H}^\perp$ is invariant for $W$ and 
 \[P_\mathcal{H} W_p \bigg|_{\mathcal{H}}=T(p), \qquad p\in P;\]
\item\label{thm.dilation.3} for every finite set $F\subset P$, 
\[
Z(F):=\sum_{U\subseteq F} (-1)^{|U|} TT^*(\vee U) \geq 0;
\]
\end{enumerate}

Here, $TT^*(\vee U)=T(s_U)T(s_U)^*$ for any $s_U\in \vee U$. By convention, $TT^*(\vee U)=0$ if $\vee U=\emptyset$.
\end{theorem}

\begin{theorem}\label{thm:diln} 
Let  $P_1$ be a submonoid of the  right-LCM monoid $P$
and assume that the inclusion $P_1\subset P$ is closed under factorization and preserves  orthogonality.
For every Nica-covariant representation $V:P_1\to \bh{H}$, there exists a Hilbert space $\cK\supset \cH$ and a Nica-covariant representation $W:P\to\bh{K}$, such that 
\begin{enumerate}
\item $\cH^\perp$ is invariant for $\{W_p: p\in P\}$;

\smallskip\item $\displaystyle P_\cH W_p\bigg|_\cH = \begin{cases} V_p, \mbox{ if } p\in P_1 \\ 0, \mbox{ if } p\notin P_1;\end{cases}$

\smallskip\item $\cH$ is reducing for $\{W_p: p\in P_1\}$, 
i.e.  $\displaystyle W_p=\begin{bmatrix}
V_p & 0 \\ 0 & D_p
\end{bmatrix}$ for  $p\in P_1$ and some $D_p$ in $\bh{\cH^\perp}$. 
\end{enumerate}
\end{theorem} 

\begin{proof} First define  $T:P\to\bh{H}$ by 
\[T(p) = \begin{cases} V_p, \mbox{ if } p\in P_1 \\ 0, \mbox{ if } p\notin P_1.\end{cases}\]
 We claim that this is a contractive representation. Clearly $\|T_p\| \leq 1$. To see $T$ is multiplicative suppose that $p,q\in P$. If any one of $p,q$ is not in $P_1$, then the factorization property states that $pq\notin P_1$ and thus $T(pq)=0=T(p)T(q)$; otherwise both $p,q\in P_1$, in which case $T(pq)=V(pq)=V_pV_q=T(p)T(q)$. 

Next we verify that the contractive representation $T$ of $P$ on $\cH$ satisfies  condition \pref{thm.dilation.3} in Theorem \ref{thm.dilation}. Let $F=\{p_1,\cdots,p_n\}$ be a finite subset of $P$ and assume $F$ contains at least one element not in $P_1$, say, $p_n\notin P_1$. Then for every  $U\subset F$ with $p_n\in U$, either $\vee U=\emptyset$ in which case $TT^*(\vee U)=0$, or else $\vee U\neq\emptyset$, in which case, we can pick a least common multiple $s_U\in \vee U \subset p_nP$ so we may  write $s_U=p_n\cdot q_{n,U}$ for some $q_{n,U}\in P$. Since $P_1$ is closed under factorization in $P$, $s_U \notin P_1$ and thus $TT^*(\vee U)=T(s_U)T(s_U)^*=0$. 
This shows that the only nontrivial terms in the sum defining $Z(F)$ in part (3) come from subsets of $F_1 := F\cap P_1$, and thus $Z(F)=Z(F_1)$. Recall that $T|_{P_1}=V$ is an isometric Nica-covariant representation, so that for all $U\subset P_1$, 
\[ \prod_{p\in U} TT^*(p) = TT^*(\vee U). \]
Then
\[
Z(F)= Z(F_1)=\sum_{U\subseteq F_1} (-1)^{|U|} \prod_{p\in U} TT^*(p) = \prod_{p\in F_1} (I-TT^*(p))\geq 0,
\]
by a straightforward expansion of the product on the right. 
This shows that $T$ satisfies part \pref{thm.dilation.3} in Theorem \ref{thm.dilation}, so by part \pref{thm.dilation.2} of the same theorem
there exists a minimal Nica-covariant isometric dilation $W$ of $T$. 

In order to see that the subspace $\cH \subset \mathcal K$ is reducing for the restriction of $W$ to $P_1$, let $p\in P_1$ 
  and
 write $W_p$  as a $2\times 2$ operator matrix with respect to $\cK=\cH\oplus \cH^\perp$: 
\[W_p=\begin{bmatrix}
T_p & 0 \\ C & D 
\end{bmatrix}.\]
Since $W_p$ is an isometry, 
\[W_p^* W_p=\begin{bmatrix}
T_p^* T_p + C^* C& C^* D\\ D^*C & D^*D
\end{bmatrix} = \begin{bmatrix}
I_{\cH} & 0 \\ 0 & I_{\cH^\perp} 
\end{bmatrix}.\]
If $p\in P_1$, then  $T_p = V_p$ is an isometry too, so $C = 0$ and $\cH$ is invariant  for $W_p$. 
\end{proof}

For inclusions that preserve orthogonality and are closed under factorization
we have the following functoriality result of the associated universal $C^*$-algebras.

\begin{corollary}\label{cor.functor.univ} Let $v$ be the universal Nica covariant representation of $P_1$ and let $w$ be that of $P$. Then the map $v_p \mapsto w_p$ for $p\in P_1$ extends to an isometric $*$-homomorphism of $C^*(P_1)$ onto the subalgebra $C^*(\{w_p: p\in P_1\})$ of $ C^*(P)$. 
\end{corollary}

\begin{proof}  By Lemma \ref{lm.respectLCM} and \cite[Theorem 3.3]{BLS2018}, there is a $*$-homomorphism $\varphi: C^*(P_1) \to  C^*(P)$ extending the given map. 
Let $v:P_1\to\bh{H}$ be a universal Nica-covariant representation. By Theorem \ref{thm:diln} there is an extension-dilation of $v$ 
 to a Nica-covariant representation $W: P\to\bh{K}$ such that $v$ is obtained from $W$ by first restricting to $P_1$ and then compressing to the reducing subspace $\cH$.
Hence
\[w_pw_q^*=\begin{bmatrix} v_pv_q^* & 0 \\ 0 & * \end{bmatrix} \qquad p,q\in P_1;\]
so the associated C*-algebra representations $\pi_W: C^*(P) \to \bh\cK$ and $\iota_v:C^*( P_1) \to \bh\cH$ satisfy
\[\pi_W(\varphi(X))=\begin{bmatrix} \iota_{v}(X) & 0 \\ 0 & * \end{bmatrix} \qquad X\in C^*(P_1).\]
Thus $\|\iota_{v}(X)\| \leq \|\pi_{W}(\varphi(X))\| \leq \|\pi_{w}(\varphi(X))\|$,
which proves that $\varphi$ is injective, hence isometric. 
\end{proof} 

\begin{theorem}\label{thm.amenable} Let $P$ be a  cancellative right-LCM monoid and suppose the inclusion $P_1 \subset P$  
is closed under factorization and preserves orthogonality. If $P$ is Nica amenable, then so is $P_1$. 
\end{theorem}

\begin{proof} 
 The map $\varphi: C^*(P_1)\to C^*(P)$ induced by the inclusion $P_1 \subset P$  is an isomorphism onto its image by Corollary \ref{cor.functor.univ}.
 Let $\lambda  : C^*(P)\to C_\lambda^*(P)$ be the left regular representation, which is faithful because we are assuming 
 $P$ to be amenable; then
the composition $\pi:= \lambda \circ\varphi : C^*(P_1)\to C_\lambda^*(P)$  is isometric and has range 
$\overline{\lspan}\{\lambda(p)\lambda(q)^*\in\mathcal{B}(\ell^2(P)): p,q\in P_1\}$.

Consider now the following commutative diagram in which the unlabeled horizontal arrows are simply inclusions.
\begin{figure}[h]
    \centering

    \begin{tikzpicture}[scale=0.9]

    \node at (-4,2) {$C^*(P_1)$};
    \node at (4,2) {$C^*(P)$};
    \node at (-4,0) {$\overline{\lspan}\{\lambda(p)\lambda(q)^*: p,q\in P_1\}$};
    \node at (4,0) {$C^*_\lambda(P)$};
    \node at (-4,-2) {$D_{P_1} \cong \overline{\lspan}\{\lambda(p)\lambda(p)^*: p,q\in P_1\} $};
    \node at (4,-2) {$D_{P,\lambda}$};
    
    \draw[->] (-2.6,2) -- (2.6,2);
    \draw[->] (-1.2,0) -- (2.6,0);
	\draw[->] (-4,1.5) -- (-4,0.5);
    \draw[->] (4,1.5) -- (4,0.5);
	\draw[->] (-.6,-2) -- (2.6,-2);
	\draw[->] (-4,-0.5) -- (-4,-1.5);
	\draw[->] (4,-0.5) -- (4,-1.5);
	
	\node at (0, 2.3) {$\varphi$};
	\node at (0, 2.3) {$\varphi$};
	\node at (0, 2.3) {$\varphi$};
	\node at (4.3, 1) {$\lambda$};
	\node at (-4.3, 1) {$\pi$};
	
	\node at (4.5, -1) {$\Phi_{P,\lambda}$};
	\node at (-4.3, -1) {$\eta$};
    \end{tikzpicture}
\end{figure}

By Proposition \ref{pro:conditionalexpectations}(1),  
$\Phi_{P,\lambda}: C^*_\lambda(P)\to D_{P,\lambda}$  is a faithful conditional expectation, and so is its restriction $\eta$
to the range of $\pi$. 
Therefore, the composition $\Phi_1=\eta\circ \pi: C^*(P_1)\to D_{P_1}$ is faithful, so $P_1$ is Nica amenable by Proposition \ref{pro:conditionalexpectations} (3). 
\end{proof} 

In the case  of Nica amenable right-LCM monoids,  functoriality holds for the reduced $C^*$-algebra.

\begin{corollary} Let $\lambda_1, \lambda$ be the left-regular representations of $P_1, P$ on $\ell^2(P_1)$ and $\ell^2(P)$, respectively. If either $P$ or $P_1$ is Nica amenable, then there is an isometric $*$-homomorphism $\psi: C^*_\lambda(P_1) \to C^*_\lambda(P)$ such that $\psi(\lambda_1(p)) = \lambda(p)$ for all $p\in P_1$. 
\end{corollary}

\begin{proof} When $P$ is Nica amenable, then so is $P_1$  by Theorem \ref{thm.amenable} and the canonical maps from $C^*(P_1)$ to $C^*_\lambda(P_1)$ and from $C^*(P)$ to $C^*_\lambda(P)$ are both isometric $*$-isomorphisms. Corollary \ref{cor.functor.univ} proves that the canonical map $\varphi:C^*(P_1)\to C^*(P)$ is an isometric $*$-homomorphism. Hence, $\phi$ must be an isometric $*$-homomorphism. 

When $P_1$ is Nica amenable,  the left-regular representation $\lambda_1:C^*(P_1)\to C^*_\lambda(P_1)$ is an isomorphism, and we consider the diagram

\begin{figure}[h]
    \centering

    \begin{tikzpicture}[scale=0.9]

    \node at (-4,2) {$C^*(P_1)$};
    \node at (4,2) {$C^*(P)$};
    \node at (-4,0) {$C^*_\lambda(P_1)$};
    \node at (4,0) {$C^*_\lambda(P)$};

    \draw[->] (-2.6,2) -- (2.6,2);
    \draw[->] (-2.6,0) -- (2.6,0);
	\draw[->] (-4.2,1.5) -- (-4.2,0.5);
    \draw[->] (4,1.5) -- (4,0.5);

	\node at (0, 2.3) {$\varphi$};
	\node at (4.3, 1) {$\lambda$};
	\node at (-4.5, 1) {$\lambda_1$};
	\node at (0, 0.3) {$\phi$};
	
	\node at (6,1) {};

    \end{tikzpicture}
\end{figure}

Corollary \ref{cor.functor.univ} proves that $\varphi: C^*(P_1)\to C^*(P)$ is an isometric $*$-homomorphism. The left-regular representation $\lambda:C^*(P)\to C^*_\lambda(P)$ is always a surjective $*$-homomorphism. Therefore, we obtain a $*$-homomorphism $\phi: C^*_\lambda(P_1)\to C^*_\lambda(P)$ by sending $\lambda_1(p)\in C^*_\lambda(P_1)$ to $\lambda(p)\in C^*_\lambda(P)$. 

Since $P_1$ is closed under factorization in $P$, $\ell^2(P_1)\subseteq \ell^2(P)$ is a reducing subspace for $\phi(C^*_\lambda(P_1))$. Therefore, the map $\phi$ must be isometric. 
\end{proof}


\section{Applications}
\label{sec:applications}
\subsection{Artin Monoids}
An  $n\times n$ Coxeter matrix $M$ is a symmetric matrix with $m_{i,i}=1$ and $m_{i,j}\in\{2,3,\cdots,\infty\}$. Denote by $\langle st\rangle^m $ the alternating product $sts...$ with a total of $m$ letters, starting with $s$ and ending with $t$ if $m$ is even, and with $s$ if $m$ is odd.
The Artin group $A_M$ associated with a Coxeter matrix  $M$ is, by definition,  the group with presentation 
\[
\Big\langle s_i: i =1, 2\, \ldots n\mid \langle s_i s_j\rangle^{m_{i,j}} = \langle s_j s_i\rangle^{m_{j,i}}   \Big\rangle;
\]  
i.e. the group generated by $s_1,\cdots, s_n$ with the relations $\langle s_i s_j\rangle^{m_{i,j}} = \langle s_j s_i\rangle^{m_{j,i}}$, with the understanding that a relation of the form $\langle s_i s_j\rangle^{\infty} = \langle s_j s_i\rangle^{\infty}$ is vacuous. The Artin monoid $A_M^+$ is defined by the same presentation (i.e. same generators and same relations), viewed in the category of monoids. If we impose the extra relation $s_i^2=e$ for each $i$, we obtain the Coxeter group $C_M$. 
The  motivating examples are the braid group $B_{n+1}$ on $n+1$ strands, introduced by Artin. $B_{n+1}$ is the Artin group associated to the $n \times n$ matrix $M$ given by
$m_{i,j}= 3$ if $|i-j|=1$ and $m_{i,j}=2$ if $|i-j|>1$; its Coxeter group is the symmetric group $S_{n+1}$.

To illustrate various other possibilities, we list a few well-known examples of Artin  groups and their corresponding monoids. 
For  $m_{i,j}=2$ for all $i\neq j$ the Artin group is the free abelian group $\Z^n$, and 
for $m_{i,j}=\infty$ for all $i\neq j$, it is the free group $\mathbb{F}_n$.  If
$M=\begin{bmatrix}1 & n \\ n & 1\end{bmatrix}$ with  $2\leq n < \infty$; the corresponding Artin group $I_n$ has Coxeter group equal to the dihedral group $D_{2n}$. More generally, 
an Artin group is of {\em spherical type} (often called of finite type, although this terminology is not universally accepted) if the corresponding Coxeter group is finite; these are  characterized by the familiar Coxeter-Dynkin diagrams, see e.g. \cite{GroveBensonBook}.  For instance $I_n$ and $B_{n+1}$ are spherical Artin groups.
The {\em right-angled} Artin groups and monoids,  $A_M$ and $A_M^+$, respectively,  
arise from matrices $M$ with $m_{i,j}\in\{2,\infty\}$. 
They are also called  graph-groups and graph-monoids, and are the  graph products 
$\Gamma\mathbb{N}$ of copies of $\Z$ and of $\N$, respectively, over the simplicial 
graph $\Gamma$ determined by having an edge $\{i,j\}$ whenever $m_{i,j} =2$.

All Artin monoids embed in their groups \cite{Paris2002}. The spherical and right-angled Artin monoids induce quasi-lattice orders in their respective groups see \cite{BrieskornSaito1972} and \cite{CrispLaca2002}, respectively, but 
 it is not known whether the other Artin monoids 
 induce quasi-lattice orders in their groups. Nevertheless, they are cancellative right LCM monoids,  in fact, weak quasi lattice orders \cite[Lemma 2.1]{BrieskornSaito1972}.
The following application of our main result completes the classification of Artin monoids with respect to Nica amenability.

\begin{theorem}\label{thm.Artin} An Artin monoid $A_M^+$ is Nica amenable if and only if it is right-angled.
\end{theorem} 

\begin{proof} Right-angled Artin monoids are Nica amenable by \cite[Theorem 20]{CrispLaca2002}. Suppose $A_M^+$ is not right-angled, then there exists $m_{i,j}\notin\{2,\infty\}$. Without loss of generality, assume $m_{1,2} = n\notin\{2,\infty\}$. Let $M_1=\begin{bmatrix}
1 & n \\ n & 1
\end{bmatrix}$ and consider the corresponding Artin monoid $A_{M_1}^+$, which embeds in $A_M^+$ by \cite[Theorem 1.3]{Crisp1999}. It is clear that $A_{M_1}^+$ is a submonoid of $A_M^+$ that is closed under factorization and respects orthogonality because $s \vee t = \left<st\right>^{m_{s\,t}}$ by \cite[Lemma 2.1]{BrieskornSaito1972}. Since $A_{M_1}^+$ is a finite-type Artin monoid that is not abelian,  $A_{M_1}^+$ is not Nica amenable by \cite[Proposition 28]{CrispLaca2002}, and by Theorem \ref{thm.amenable}, this implies that $A_M^+$ is  not Nica amenable .
\end{proof} 


\subsection{Graph Products}
One way to construct new right-LCM monoids is to take  graph products of old ones. Let $\Gamma$ be a simple graph on $n$ vertices, and for each vertex $1\leq i\leq n$, let $P_i$ be a right-LCM monoid. We construct the graph product $\Gamma P_i$ by taking the free product monoid modulo the relations that make  $P_i$ commute with $P_j$ whenever $(i,j)$ is an edge of $\Gamma$. Then $\Gamma P_i$ is a right-LCM monoid by \cite[Theorem 2.6]{FK2009}, see also \cite[Theorem 5.25]{ABCD2019}. 

\begin{example} When $\Gamma$ is a complete graph on $n$ vertices, the graph product $\Gamma P_i$ is simply the direct product of the $P_i$. 
When  $\Gamma$ contains no edges, the graph product $\Gamma P_i$ is their free product.
\end{example}

In the case when the $(G_i, P_i)$ are quasi-lattice ordered, it is known \cite[Theorem 20]{CrispLaca2002} that the graph product $(\Gamma G_i, \Gamma P_i)$ is Nica amenable if all the $G_i$ are amenable groups.
We have a partial converse to this.

\begin{corollary} \label{cor:graphprod}
Suppose $\Gamma $ is a simplicial graph and let $P_k$ be a cancellative right-LCM monoid for each vertex $k$ of $\Gamma$. If the graph product $\Gamma_i P_i$ is Nica amenable, then all the $P_i$ are Nica amenable.
\end{corollary}

\begin{proof}  It is clear that each $P_k$ embeds in $\Gamma P_i$ as a sub-monoid that is closed under factorization. By \cite[Lemma 2.7]{FK2009}  each $P_k$ preserves  orthogonality. The result then follows from Theorem~\ref{thm.amenable}. 
\end{proof} 

\begin{remark} Whether the amenability of each $P_i$ is sufficient for the amenability of their graph product remains an open problem. The known results in this direction require  $P_i$ to embed  in an amenable group $G_i$ for each $i$ \cite{ABCD2019, CrispLaca2002}. 
\end{remark}



\begin{thebibliography}{10}

\bibitem{ABCD2019}
A.~an~Huef, B.~Nucinkis, C. F.~Sehnem, and D.~Yang.
\newblock Nuclearity of semigroup {C}*-algebras.
\newblock {\em https://arxiv.org/abs/1910.04898}, 2019.

\bibitem{HRT2018}
A.~an~Huef, I.~Raeburn, and I.~Tolich.
\newblock H{NN} extensions of quasi-lattice ordered groups and their operator
  algebras.
\newblock {\em Doc. Math.}, 23:327--351, 2018.

\bibitem{BrieskornSaito1972}
E.~Brieskorn and K.~Saito.
\newblock Artin-{G}ruppen und {C}oxeter-{G}ruppen.
\newblock {\em Invent. Math.}, 17:245--271, 1972.

\bibitem{BLS2017}
N.~Brownlowe, N.~S. Larsen, and N.~Stammeier.
\newblock On {$C^*$}-algebras associated to right {LCM} semigroups.
\newblock {\em Trans. Amer. Math. Soc.}, 369(1):31--68, 2017.

\bibitem{BLS2018}
N.~Brownlowe, N.~S. Larsen, and N.~Stammeier.
\newblock {$C^*$}-algebras of algebraic dynamical systems and right {LCM}
  semigroups.
\newblock {\em Indiana Univ. Math. J.}, 67(6):2453--2486, 2018.

\bibitem{Coburn1967} 
L. A. Coburn.
\newblock The $C^*$-algebra generated by an isometry I.
\newblock  {\em Bull. Amer. Math. Soc.} {73}:722--726,  1967. 
 
\bibitem{Crisp1999}
J.~Crisp.
\newblock Injective maps between {A}rtin groups.
\newblock In {\em Geometric group theory down under ({C}anberra, 1996)}, pages
  119--137. de Gruyter, Berlin, 1999.

\bibitem{CrispLaca2002}
J.~Crisp and M.~Laca.
\newblock On the {T}oeplitz algebras of right-angled and finite-type {A}rtin
  groups.
\newblock {\em J. Aust. Math. Soc.}, 72(2):223--245, 2002.

\bibitem{CrispLaca2007}
J.~Crisp and M.~Laca.
\newblock Boundary quotients and ideals of {T}oeplitz {$C^*$}-algebras of
  {A}rtin groups.
\newblock {\em J. Funct. Anal.}, 242(1):127--156, 2007.

\bibitem{Cuntz1977} J. Cuntz, 
\newblock Simple $C^*$-algebras generated by isometries.
 \newblock {\em   Comm. Math. Phys.}  {57}:173--185, 1977. 


\bibitem{Cuntz1981} 
J. Cuntz, 
\newblock {K-theory for certain C*-algebras}.
\newblock {\em Annals of Mathematics} {113}:181--197, 1981.

\bibitem{CDL2013}
J.~Cuntz, C.~Deninger, and M.~Laca.
\newblock {$C^*$}-algebras of {T}oeplitz type associated with algebraic number
  fields.
\newblock {\em Math. Ann.}, 355(4):1383--1423, 2013.

\bibitem{CELY2017}
J. Cuntz, S. Echterhoff, X.~Li, and G. Yu.
\newblock K-Theory for group C*-algebras and semigroup C*-algebras.
\newblock Oberwolfach Seminars Vol. 47, Birkh\"auser/Springer 2017.

\bibitem{Day1957}
M.M. Day.
\newblock Amenable semigroups.
\newblock {\em Illinois J. Math.} 1(4):509-544, 1957.

\bibitem{Douglas1972} 
R.G. Douglas.
\newblock On the $C^*$-algebra of a one-parameter semigroup of isometries
 \newblock {\em Acta Math.} {128}:143--152, 1972. 
  
\bibitem{ELR2016} 
S. Eilers,  X. Li,  and  E. Ruiz. 
\newblock The isomorphism problem for semigroup C?-algebras of right-angled Artin monoids.
\newblock  {\em Doc. Math.}  21  (2016), 309--343.
 
\bibitem{Exel1997}
R.~Exel.
\newblock Amenability for {F}ell bundles.
\newblock {\em J. Reine Angew. Math.}, 492:41--73, 1997.

\bibitem{FK2009}
J.~Fountain and M.~Kambites.
\newblock Graph products of right cancellative monoids.
\newblock {\em J. Aust. Math. Soc.}, 87(2):227--252, 2009.

\bibitem{GroveBensonBook}
L.~C. Grove and C.~T. Benson.
\newblock {\em Finite reflection groups}, volume~99 of {\em Graduate Texts in
  Mathematics}.
\newblock Springer-Verlag, New York, second edition, 1985.

\bibitem{Ivanov2010}
N. Ivanov.
\newblock {\em The K-theory of Toeplitz C?-algebras of right-angled Artin groups} 
\newblock Trans. Amer. Math. Soc. 362 (2010), no. 11, 6003--6027.

\bibitem{LacaRaeburn1996}
M.~Laca and I.~Raeburn.
\newblock Semigroup crossed products and the {T}oeplitz algebras of nonabelian
  groups.
\newblock {\em J. Funct. Anal.}, 139(2):415--440, 1996.

\bibitem{BLi2019}
B.~Li.
\newblock Regular dilation and {N}ica-covariant representation on right {LCM}
  semigroups.
\newblock {\em Integr. Equ. Oper. Theory}, 91(4):36, 2019.

\bibitem{XLi2012}
X.~Li.
\newblock Semigroup {${C}^*$}-algebras and amenability of semigroups.
\newblock {\em J. Funct. Anal.}, 262(10):4302--4340, 2012.


\bibitem{Nica1992}
A.~Nica.
\newblock {$C\sp *$}-algebras generated by isometries and {W}iener-{H}opf
  operators.
\newblock {\em J. Operator Theory}, 27(1):17--52, 1992.

\bibitem{Norling2014} 
M.~D.~Norling.
\newblock Inverse semigroup {${C}^*$}-algebras associated with left cancellative semigroups.
\newblock {\em Proc. Edinb. Math. Soc.}, 57:533--564, 2014. 

\bibitem{LOS2019}
X.~Li, T. Omland, and J. Spielberg.
\newblock C*-algebras of right LCM one-relator monoids and Artin-Tits groups of finite type.
\newblock preprint arXiv:1807.08288.

\bibitem{Paris2002} 
L. Paris
\newblock Artin monoids inject in their groups.
\newblock {\em Comment. Math. Helv.}, 77:609--637, 2002.

\bibitem{Popescu1989}
G.~Popescu.
\newblock Isometric dilations for infinite sequences of noncommuting operators.
\newblock {\em Trans. Amer. Math. Soc.}, 316(2):523--536, 1989.


\bibitem{Starling2015}
C.~Starling.
\newblock Boundary quotients of {$\rm C^*$}-algebras of right {LCM} semigroups.
\newblock {\em J. Funct. Anal.}, 268(11):3326--3356, 2015.

\end{thebibliography}
\bibliographystyle{amsalpha}


\end{document}